\documentclass[12pt,a4paper]{amsart}
\usepackage{geometry}
\usepackage{tabls}

\usepackage{url}

\usepackage{amssymb}
\usepackage{mathrsfs}
\usepackage{amsmath}
\usepackage{amsthm}
\usepackage{amscd}
\usepackage{fancyhdr}
\usepackage{enumerate}
\usepackage{paralist}
\usepackage{graphicx}
\usepackage{cite}

\usepackage{footmisc}
\numberwithin{equation}{section}

\theoremstyle{plain}

\newtheorem{Thm}[equation]{Theorem}
\newtheorem{lem}[equation]{Lemma}

\newtheorem{rem}[equation]{Remark}
\newtheorem{ex}[equation]{Example}

\begin{document}

\title{Rational deformations of the set of multiple zeta-star values}

\author{Jiangtao Li}

\email{lijiangtao@csu.edu.cn}
\address{Jiangtao Li \\ School of Mathematics and Statistics, HNP-LAMA, Central South University, Hunan Province, China}

\begin{abstract}
   In this paper we study the derived sets for the rational deformations of multiple zeta-star values. By using the theory of bounded variation functions, we will give  function decompositions which describe  the metric structure of the derived sets. The connection between the rational deformation of multiple zeta-star values and the $n$-th Cantor set in fractal geometry is also discussed.
    
    \end{abstract}

\let\thefootnote\relax\footnotetext{
2020 $\mathnormal{Mathematics} \;\mathnormal{Subject}\;\mathnormal{Classification}$. 11K55, 11M32. \\
$\mathnormal{Keywords:}$ Hausdorff dimension, Multiple zeta-star values.
 }

\maketitle

\section{Introduction}\label{int}
      Multiple zeta values are defined by
      \[
      \zeta(k_1,\cdots,k_r)=\sum_{n_1>\cdots>n_r>0}\frac{1}{n_1^{k_1}\cdots n_r^{k_r}}, k_1\geq 2,k_2, \cdots, k_{r}\geq 1.
      \]
      For $k_1\geq 2,k_2, \cdots, k_{r}\geq 1$, the multiple zeta-star value $\zeta^{\star}(k_1,\cdots,k_r)$ is defined by
        \[
        \zeta^{\star}(k_1,\cdots,k_r)=\sum_{n_1\geq\cdots\geq n_r\geq 1}\frac{1}{n_1^{k_1}\cdots n_r^{k_r}}.
        \]
      For multiple zeta value $\zeta(k_1,\cdots,k_r)$ and multiple zeta-star value $\zeta^\star(k_1,\cdots,k_r)$, denote by $N=k_1+\cdots+k_r$ and $r$ its weight  and depth respectively. 
      
       For the algebra structure of multiple zeta values, now there are many related results  \cite{bro}, \cite{eca}, \cite{hof},\cite{ikz},\cite{ow},\cite{oz}. On the other hand, for the topology and metric structure of multiple zeta value and multiple zeta-star values we still known very little.
      
      In \cite{kum}, Kumar investigated the order structure and the $n$-th derived set of the set of multiple zeta values. In \cite{liz}, Li and Pan obtained the derived sets for some $q$-analogue of multiple zeta values. In \cite{GLW}, Gorodetsky, Lichtman and Wong   calculated a limit of a sequence of multiple integrals.   In \cite{lit}, the author found that there is a natural total order structure on the set of multiple zeta-star values. Furthermore, the author showed that:  Define 
 \[
 \mathcal{T}=\Big{\{} (k_1,k_2, \cdots, k_r,\cdots)\,\Big{|}\,k_1\geq 2,k_i\geq 1, i\geq 2, k_s\geq 2 \mathrm{\;for \; some} \; s\geq 2 \;\mathrm{if} \;k_1=2    \Big{\}}.
 \]
 Then the map 
 \[
 \eta: \mathcal{T}\rightarrow (1, +\infty),
  \]
    \[
  {\bf k}=(k_1,k_2, \cdots, k_r,\cdots) \mapsto x=\mathop{\mathrm{lim}}_{r\rightarrow +\infty}\zeta^{\star}(k_1,k_2,\cdots,k_r)  \]
 is bijective.   As a result, $\left(  \mathcal{Z}^{\star}   \right)^{(n)}= [1, +\infty], \forall\, n\geq 1$.   Furthermore, by using the map $\eta$, the author showed that every mutiple zeta-star value is not integer although the set of multiple zeta-star values is dense in $(1,+\infty)$. In \cite{hmo}, M. Hirose, H. Murahara and T. Onozuka also found the bijective map $\eta$ independently. They also calculate the images of many special sequences of non-negative integers.
 
 Many basic properties of the map $\eta$ are discussed in \cite{lit}.  But it seems that the map $\eta$ is still quite mysterious.    As all the multiple zeta-star values are conjectured to be transcendental numbers, it is hard to understand the metric structure of the map $\eta$ directly. 
 
  For $k_1\geq 2, k_2,\cdots, k_r\geq 1$, the multiple zeta-star value $\zeta^{\star}(k_1,k_2,\cdots,k_r)  $ has integral representation of the form 
 \[
 \begin{split}
 &\;\;\;\;\zeta^{\star}(k_1,k_2,\cdots,k_r)\\
 &=\mathop{\int}_{[0,1]^{k_1+\cdots+k_r}} \frac{dx_1dx_2\cdots dx_{k_1+\cdots+k_r}}{(1-x_1\cdots x_{k_1})(1-x_1\cdots x_{k_1+k_2}) \cdots (1-x_1\cdots x_{k_1+k_2+\cdots +k_r}) }.\\
 \end{split}
  \]
  Define 
  \[
 \begin{split}
 &\;\;\;\;T_n(k_1,k_2,\cdots,k_r)\\
 &=\mathop{\int}_{[0,1]^{k_1+\cdots+k_r}}\left[1+x_1\cdots x_{k_1}+ \cdots +(x_1\cdots x_{k_1})^n\right]\left[1+x_1\cdots x_{k_1+k_2}+\cdots+ (x_1\cdots x_{k_1+k_2})^n\right] \\
 &\;\;\;\; \cdots \left[ 1+x_1\cdots x_{k_1+k_2+\cdots +k_r}+\cdots+( x_1\cdots x_{k_1+k_2+\cdots +k_r})^n\right] dx_1dx_2\cdots dx_{k_1+\cdots+k_r}.\\ \end{split}
  \]
  It is clear that $T_n(k_1,k_2,\cdots,k_r)$ is rational and 
  \[
  \zeta^{\star}(k_1,k_2,\cdots,k_r)=\mathop{\lim}_{n\rightarrow +\infty} T_n(k_1,k_2,\cdots,k_r) .   \]  
  Thus the number $T_n(k_1,k_2,\cdots,k_r)$ can be viewed as  a rational deformation of the multiple zeta-star value $\zeta^{\star}(k_1,k_2,\cdots,k_r)$.
  
  In this paper, in order to give a glimpse of the metric structure of the map $\eta$, we study the derived sets for the rational deformations of the set of multiple zeta-star values. Define 
  \[
 \widehat{\mathcal{T}}=\Big{\{} (k_1,k_2, \cdots, k_r,\cdots)\,\Big{|}\,k_i\geq 1, i\geq 1,k_i\in\mathbb{N}   \Big{\}},
 \]
  \[
 \widehat{\mathcal{T}}_p=\Big{\{} (k_1,k_2, \cdots, k_r,\cdots)\,\Big{|}\,  k_i\geq p, i\geq 1,k_i\in\mathbb{N}   \Big{\}},
 \]
 \[
 \widehat{\mathcal{T}}^{R}=\Big{\{} (k_1,k_2, \cdots, k_r,\cdots)\,\Big{|}\,k_i\geq 1, i\geq 1,k_s\geq 2, k_i=1\; \mathrm{if} \; i>s\; \mathrm{for \;some} \;s  \Big{\}}.
 \]   For $k_1, k_2\in \mathbb{N}$, define a total order $\succ$ on $\mathbb{N}$ such that $k_1\succ k_2$ iff $k_1<k_2$. We equip $\widehat{\mathcal{T}}$ with the lexicographical order $\succ$ which is induced by the order $\succ$ on $\mathbb{N}$.  By the binary expression of real numbers, it is clear that the map
 \[
 \tau:  \widehat{\mathcal{T}}\rightarrow (0, 1], \]
 \[
 (k_1,k_2,\cdots,k_r,\cdots)\mapsto \frac{1}{2^{k_1}}+\frac{1}{2^{k_1+k_2}}+\cdots+\frac{1}{2^{k_1+\cdots+k_r}}+\cdots
 \] 
 is bijective and $\tau({\bf k})> \tau({\bf l})$ if and only if    ${\bf k}\succ {\bf l}$. It is easy to see that
 \[
 \tau(\widehat{\mathcal{T}}^{R} )=\big{\{}x\in (0,1)\,\big{|}\, 2^l x\in \mathbb{N} \;\mathrm{for \;some} \;l \in \mathbb{N}\big{\}} .\]
 The set $\tau(\widehat{\mathcal{T}}^{R} )$ is a countable dense subset of $(0, 1)$. By Lemma $3.4$ in \cite{lit}, it is clear that
 \[
\Big{\{} \mathop{\lim}_{n\rightarrow+\infty}\eta_n\left(  {\bf k} \right)\,\Big{|}\, {\bf k}\in\widehat{\mathcal{T}}^{R} \cap \mathcal{T} \Big{ \}}=\mathcal{Z}^\star,
 \]
 where $\mathcal{Z}^\star$ is the set of multiple zeta-star values.
   \begin{Thm}\label{rat} Define the map \[\eta_n: \widehat{\mathcal{T}}\rightarrow (1, +\infty)\]
 as \[
 (k_1,k_2, \cdots, k_r,\cdots)\mapsto \mathop{\lim}_{r\rightarrow +\infty} T_n(k_1,k_2,\cdots,k_r). \]
 $(i)$ For $n\geq 1$, $$\eta_n(\widehat{\mathcal{T}}- \widehat{\mathcal{T}}^{R})=(1,e^{1+\frac{1}{2}+\cdots+\frac{1}{n}}], \eta_n(\widehat{\mathcal{T}}^{R})\subseteq(1,e^{1+\frac{1}{2}+\cdots+\frac{1}{n}})$$ and $\eta_n\left((1,1,\cdots,1, \cdots) \right)=e^{1+\frac{1}{2}+\cdots+\frac{1}{n}}$;\\
  $(ii)$ For ${\bf k,l}\in  \widehat{\mathcal{T}}_2$, we have $\eta_n({\bf k})> \eta_n({\bf l})$ if and only if    ${\bf k}\succ {\bf l}$.\\
   $(iii)$ Denote by $\mathrm{dim}_{\mathrm{H}} C$ the Hausdorff dimension of $C$ for any $C\subseteq \mathbb{R}$, then 
  \[
  \mathrm{dim}_{\mathrm{H}} \eta_n(  \widehat{\mathcal{T}}_2)=\frac{\mathrm{log}\; \frac{1+\sqrt{5}}{2}}{\mathrm{log}\;2} .\] \end{Thm}
  Although the natural logarithm $e$ does not appear in the theory of multiple zeta values, it appears naturally in the derived set of rational deformations of multiple zeta-star values. Beware that the statement $(ii)$ does not hold in general cases.
  
 For convenience we identify $\widehat{\mathcal{T}}$ with $(0, 1]$ by the map $\tau$. For the map $$F=\eta\circ \tau^{-1}|_{(0,\frac{1}{2})}: (0,\frac{1}{2})\rightarrow (1,+\infty),$$ from the results in \cite{lit} we have that $F$ is a topology homemorphism.
 For the map $\eta_n$, one has the following  decomposition.
   \begin{Thm}\label{dec}
  Define $F_n: (0, 1]\rightarrow (1, +\infty)$ as $F_n=\eta_n\circ \tau^{-1}$. Then $F_n$ is a bounded variation function and $$F_n=G_n-H_n.$$ Here $$G_n,H_n:(0,1]\rightarrow (1,+\infty)$$ are strictly increasing functions. Furthermore, $G_n$ is continuous, $H_n$ is discontinuous at $\tau(\widehat{\mathcal{T}}^{R} )$ and continuous at $\tau(\widehat{\mathcal{T}} -\widehat{\mathcal{T}}^{R} )$. What is more, the image $\mathrm{Im}\,(H_n)$ is bi-Lipschitz      equivalent to    $$C_{\frac{1}{n+2}}-  C_{\frac{1}{n+2}}^b.$$ Here $C_{\frac{1}{n+2}}$ is the $(n+2)$-th Cantor set and $C_{\frac{1}{n+2}}^b$ is a countable subset of $C_{\frac{1}{n+2}}$ .  As a result,     \[  \mathrm{dim}_{\mathrm{H}} \mathrm{Im}\,(H_n)=\frac{\mathrm{log}\,2}{\mathrm{log}(n+2)} .\] \end{Thm}
 
  \section{Hausdorff dimension and Cantor set}
  In this section, we will review the definitions of Hausdorff dimension and Cantor set. The reference is \cite{fal}.

 For $F\subseteq \mathbb{R}^m$, let
 \[
 \mathcal{H}_{\infty}^s(F)=\mathrm{inf}\Big{\{} \sum_{i=1}^{\infty} |U_i|^s\,\big{|}\, \{U_i\} \,\mathrm{is\; a\; countable\;cover\; of}\;F  \Big{\}}.
 \]
 Define the Hausdorff dimension of $F$ as 
 \[
 \mathrm{dim}_{\mathrm{H}}F=\mathrm{inf}\{s\geq 0: \mathcal{H}^s_{\infty}(F)=0\}.
 \]
 \begin{ex}
 If $F$ is a countable subset of $\mathbb{R}^m$, then  $\mathrm{dim}_{\mathrm{H}}F=0$.
  \end{ex}
  \begin{ex}
  If $F=[0,1]^m\subseteq \mathbb{R}^m$, then $\mathrm{dim}_{\mathrm{H}}F=m$.  \end{ex}
  
  For $m\geq 3$, the $m$-th Cantor set $C_{\frac{1}{m}}$ is defined as 
  \[
  C_{\frac{1}{m}}=\Big{\{} x\in [0,1]\,\Big{|}\, x=\frac{a_1}{m}+\frac{a_2}{m^2}+\cdots+\frac{a_r}{m^r}+\cdots, a_i\in \{0, m-1\}  \Big{\}}.
  \]
  One can check that $ C_{\frac{1}{m}}$ is a totally disconnected closed subset of $[0,1]$. From the theory of iterated function systems (Theorem $9.3$ in \cite{fal}), it follows that
  \[
  \mathrm{dim}_{\mathrm{H}} \mathrm{Im}\,(C_{\frac{1}{m}})=\frac{\mathrm{log}\,2}{\mathrm{log}\,m} .
    \]
    Define 
     \[
  C_{\frac{1}{m}}^b=\Big{\{} x\in [0,1]\,\Big{|}\, x=\frac{a_1}{m}+\frac{a_2}{m^2}+\cdots+\frac{a_r}{m^r}, a_i\in \{0, m-1\} \; \mathrm{for} \;\mathrm{some}\; r \Big{\}}.
  \]
  
    \section{The order structure of the image}
    There is a natural order structure on the set $ \widehat{\mathcal{T}} $. In this section we will study the order structure of the image of $\eta_n$. We will show that  the order structure of the image of $\eta_n$ is not compatible with  the  order structure on the set $ \widehat{\mathcal{T}} $ in general cases. We also calculate the Hausdorff dimension for a subset of the image.   
    \begin{lem}\label{leq}
    For $r\geq 2$ and $u\in (0,+\infty)$, 
     \[\mathop{\int}_{[0,1]^{r-1}}\left[1+ux_1\right]\cdots \left[1+ux_1\cdots x_{r-1}\right] \left(x_1\cdots x_{r-1}\right)^mdx_1dx_2\cdots dx_{r-1} < \frac{e^{(m+1)u}}{(m+1)^{r-1}} .\]    \end{lem}
      \noindent{\bf Proof:}
      We have 
      \[
      \begin{split}
      &\;\;\;\;\mathop{\int}_{[0,1]^{r-1}}\left[1+ux_1\right]\cdots \left[1+ux_1\cdots x_{r-1}\right] \left(x_1\cdots x_{r-1}\right)^mdx_1dx_2\cdots dx_{r-1}   \\
    &=   \mathop{\int}_{[0,1]^{r-1}}\left[1+\sum_{s=1}^{r-1}\sum_{1\leq i_1< \cdots<i_s\leq r-1}(x_1\cdots x_{i_1})(x_1\cdots x_{i_2})\cdots (x_1\cdots x_{i_s}) u^s       \right] \\
    &\;\;\;\;\;\; \;\;\;\;\;\;\;\;\;\;\;\;\;\;\;\;\;\;\cdot \left(x_1x_2\cdots x_{r-1}  \right)^m  dx_1dx_2\cdots dx_{r-1}  \\
    \end{split}
    \]
    \[
    \begin{split}
    & = \mathop{\int}_{[0,1]^{r-1}}\Big{[}\left(x_1x_2\cdots x_{r-1}\right)^m +\sum_{s=1}^{r-1}\sum_{1\leq i_1< \cdots<i_s\leq r-1}  (x_1\cdots x_{i_1})^{m+s}(x_{i_1+1}\cdots x_{i_2})^{m+s-1}\cdots \\
   & \;\;\;\;\;\;\;\;\;\;\;\;\;\;(x_{i_{s-1}+1}\cdots x_{i_s})^{m+1} \left(x_{i_s+1}\cdots x_{r-1}\right)^m\cdot  u^s      \Big{] } dx_1dx_2\cdots dx_{r-1}  \\
   &=\left(\frac{1}{(m+1)^{r-1}}+\sum_{s=1}^{r-1}\sum_{1\leq i_1<\cdots<i_s\leq r-1}\frac{ u^s}{(m+s+1)^{i_1}\cdots (m+2)^{i_s-i_{s-1}}(m+1)^{r-1-i_s}}  \right)\\
   &= \frac{1}{(m+1)^{r-1}} \left(1+\sum_{s=1}^{r-1}\sum_{ 1\leq i_1<\cdots<i_s\leq r-1}\left(\frac{m+1}{m+s+1}\right)^{i_1}\cdots \left( \frac{m+1}{m+2} \right)^{i_s-i_{s-1}} u^s \right)\\
      &< \frac{1}{(m+1)^{r-1}} \left(1+\sum_{s=1}^{+\infty}\sum_{ 1\leq i_1<\cdots<i_s< +\infty}\left(\frac{m+1}{m+s+1}\right)^{i_1}\cdots \left( \frac{m+1}{m+2} \right)^{i_s-i_{s-1}} u^s \right)\\
      &\leq  \frac{1}{(m+1)^{r-1}} \left(1+\sum_{s=1}^{+\infty}\frac{\frac{m+1}{m+s+1}}{1-\frac{m+1}{m+s+1}}\cdots \frac{\frac{m+1}{m+2}}{1-\frac{m+1}{m+2}}
       u^s \right)\\
       &\leq \frac{e^{(m+1)u}}{(m+1)^{r-1}} . \end{split}
    \]
      $\hfill\Box$\\
      
          \begin{lem}\label{111}
    For $u\in \mathbb{R}$, define $$t_n(x)=x+x^2+\cdots+x^n,$$
    \[
    \begin{split}
    &f_n(u)=\mathop{\lim}_{r\rightarrow +\infty} \mathop{\int}_{[0,1]^{r}}\left[1+t_n(ux_1)\right]\left[1+t_n(ux_1x_{2})\right]\cdots \left[ 1+t_n(ux_1\cdots x_{r})\right] dx_1dx_2\cdots dx_{r}.
    \end{split}    \]
    Then 
    \[
    f_n(u)=\mathrm{exp}\,\left(u+\frac{u^2}{2}+\cdots+\frac{u^n}{n}\right).
    \]
    \end{lem}
    \noindent{\bf Proof:} Firstly we need to show that for $u\in (0,+\infty)$, the limit 
    \[
    \mathop{\lim}_{r\rightarrow +\infty} \mathop{\int}_{[0,1]^{r}}\left[1+t_n(ux_1\right]\left[1+t_n(ux_1x_{2})\right]\cdots \left[ 1+t_n(ux_1\cdots x_{r})\right] dx_1dx_2\cdots dx_{r}    \]
    exists. For fixed $u$ and $n$, define $M_r$ as 
    \[
    M_r=\mathop{\int}_{[0,1]^{r}}\left[1+t_n(ux_1)\right]\left[1+t_n(ux_1x_{2})\right]\cdots \left[ 1+t_n(ux_1\cdots x_{r})\right] dx_1dx_2\cdots dx_{r}  .    \]
    One has $M_{r-1}<M_{r}$ and 
    \[
    \begin{split}
    &\;\;\;\;M_{r}-M_{r-1}\\
    &=\mathop{\int}_{[0,1]^{r}}\left[1+t_n(ux_1)\right]\cdots \left[1+t_n(ux_1\cdots x_{r-1})\right] t_n(ux_1\cdots x_{r})dx_1dx_2\cdots dx_{r}  \\
    &< \mathop{\int}_{[0,1]^{r}}\left[1+nux_1\right]\cdots \left[1+nux_1\cdots x_{r-1}\right] nux_1\cdots x_{r}dx_1dx_2\cdots dx_{r}    \\
    &\leq \frac{nu}{2} \mathop{\int}_{[0,1]^{r-1}}\left[1+nux_1\right]\cdots \left[1+nux_1\cdots x_{r-1}\right] x_1\cdots x_{r-1}dx_1dx_2\cdots dx_{r-1}   \\
   &\leq \frac{nu}{2^r}e^{2nu}. \end{split}
    \]
    Here the last inequality follows from Lemma \ref{leq}.
    As a result, $M_r<M_1+\frac{nu}{2}e^{2nu}$. Thus the limit  \[
    \mathop{\lim}_{r\rightarrow +\infty} M_r \] exists.
    
    By changing of variables, one has 
    \[uf_n(u)=\mathop{\lim}_{r\rightarrow +\infty} \int^u_0dx_1\mathop{\int}_{[0,1]^{r-1}}\left[1+t_n(x_1)\right]\left[1+t_n(x_1x_{2})\right]\cdots \left[ 1+t_n(x_1\cdots x_{r})\right] dx_2\cdots dx_{r}.\]
    It follows that 
    \[
    \begin{split}
    \frac{d}{du}\left(uf_n(u)  \right)&=\mathop{\lim}_{r\rightarrow +\infty}\left[1+t_n(u)\right] \mathop{\int}_{[0,1]^{r-1}}\left[1+t_n(ux_{2})\right]\cdots \left[ 1+t_n(ux_2\cdots x_{r})\right] dx_2\cdots dx_{r}  \\
    &=\left[1+t_n(u)\right] f_n(u).    \end{split} \]
    Thus 
    \[
    \frac{f_n^\prime(u)}{f_n(u)}=\frac{t_n(u)}{u}.
    \]
    By integration on both sides, one has 
    \[
    f_n(u)=C\,\mathrm{exp}\,\left(u+\frac{u^2}{2}+\cdots+\frac{u^n}{n}     \right).
    \]
    Since $f_n(0)=1$, we have 
    \[
    f_n(u)=\mathrm{exp}\,\left(u+\frac{u^2}{2}+\cdots+\frac{u^n}{n}     \right).
    \]
$\hfill\Box$\\

For convenience, we write
\[
e_n(u)=\mathrm{exp}\,\left(u+\frac{u^2}{2}+\cdots+\frac{u^n}{n}     \right).
\]
The following observation will be used in the next section.
 \begin{rem} One has \[
 \int^1_0 \left[1+t_n(ux)\right] e_n(ux)dx=e_n(u).
 \]
 \end{rem}

\begin{lem}\label{two}
 For $u\in (0,1)$, define 
    \[
    \begin{split}
    &\;\;\;\;g_n(u)\\
    &=\mathop{\lim}_{r\rightarrow +\infty} \mathop{\int}_{[0,1]^{2r}}\left[1+t_n(ux_1x_{2})\right]\left[1+t_n(ux_1x_{2}x_3x_4)\right]\cdots \left[ 1+t_n(ux_1x_2\cdots x_{2r})\right] dx_1dx_2\cdots dx_{2r}.
    \end{split}    \]
    Then $$g_n(u)=1+\sum_{m= 1}^n\frac{2}{(m+1)(m+2)}u^m+\sum_{m>n} a_m u^m,$$
    \[
   \left[1+t_n(u)\right]g_n(u) = 1+\sum_{m\geq 1}(m+1)^2a_m u^m
    \]
    and \[
    \begin{split}
    &\;\;\;\;\int^1_0 \left[1+t_n(ux)\right]g_n(ux)dx \\
       &= 1+\sum_{m= 1}^n \frac{2}{m+2} u^m  + \sum_{m>n}(m+1) a_m u^m  .
    \end{split}
    \]
    Here 
    \[
    0<a_m\leq \frac{2}{m(m+2)}\left(  \frac{1}{m+1-n}-\frac{1}{m+1} \right),\;\;\;m>n.
    \]
\end{lem}
\noindent{\bf Proof:} Define 
\[
N_r=\mathop{\int}_{[0,1]^{2r}}\left[1+t_n(ux_1x_{2})\right]\left[1+t_n(ux_1x_{2}x_3x_4)\right]\cdots \left[ 1+t_n(ux_1x_2\cdots x_{2r})\right] dx_1dx_2\cdots dx_{2r},\]
we have $N_r<N_{r+1}$ and 
\[
N_r<\mathop{\int}_{[0,1]^{r}}\left[1+t_n(ux_{2})\right]\left[1+t_n(ux_{2}x_4)\right]\cdots \left[ 1+t_n(ux_2\cdots x_{2r})\right] dx_2dx_4\cdots dx_{2r}=M_r.\]
As a result,   the limit  \[
    \mathop{\lim}_{r\rightarrow +\infty} N_r \] exists.
By changing of variables, one has 
    \[
    \begin{split}
    &\;\;\;\;ug_n(u)\\
    &=\mathop{\lim}_{r\rightarrow +\infty} \int^u_0dx_1\mathop{\int}_{[0,1]^{2r-1}}\left[1+t_n(x_1x_2)\right]\left[1+t_n(x_1x_{2}x_3x_4)\right]\cdots \left[ 1+t_n(x_1x_2\cdots  x_{2r})\right] dx_2\cdots dx_{2r}.\\
    \end{split}
    \]
    By taking derivatives, one has 
    \[
    \frac{d}{du}\left( ug_n(u)  \right)=\mathop{\lim}_{r\rightarrow +\infty}\mathop{\int}_{[0,1]^{2r-1}}\left[1+t_n(u x_2)\right]\left[1+t_n(u x_{2}x_3x_4)\right]\cdots \left[ 1+t_n(u x_2\cdots  x_{2r})\right] dx_2\cdots dx_{2r}.    \]
    Again, by changing of variables, one has
     \[
     \begin{split}
   &\;\;\;\;u\cdot  \frac{d}{du}\left( ug_n(u)  \right)\\
   &=\mathop{\lim}_{r\rightarrow +\infty}\int^u_0dx_2\mathop{\int}_{[0,1]^{2r-2}}\left[1+t_n( x_2)\right]\left[1+t_n( x_{2}x_3x_4)\right]\cdots \left[ 1+t_n( x_2\cdots  x_{2r})\right] dx_3\cdots dx_{2r}.  \\
   \end{split}  \]
   By taking derivatives, one has 
   \[
   \begin{split}
   &\;\;\;\;\frac{d}{du}\left[u\cdot  \frac{d}{du}\left( ug_n(u)  \right)\right]\\
   &=\mathop{\lim}_{r\rightarrow +\infty}\mathop{\int}_{[0,1]^{2r-2}}\left[1+t_n( u)\right]\left[1+t_n( ux_3x_4)\right]\cdots \left[ 1+t_n( u x_3\cdots  x_{2r})\right] dx_3\cdots dx_{2r} \\
   &=\left[1+t_n( u)\right] g_n(u).   \end{split}   \]
   The above formula is equivalent to 
   \[
   u^2 g_n^{\prime\prime}(u)+3ug_n^{\prime}(u)+g_n(u)=[1+{t_n(u)}] g_n(u). \tag{1}
   \]
   Since $g_n(0)=1$, by the  formula $(1)$ we have $g_n^{\prime}(0)=\frac{1}{3}$.
   Assuming that 
   \[
   g_m(u)=a_0+\sum_{m\geq 1}a_m u^m, a_0=1, a_1=\frac{1}{3},
   \]
   by the  formula $(1)$, one has  
   \[
   \begin{split}
   &\;\;\;\;1+\sum_{m\geq 1}\left[ m(m-1)a_m+3ma_m+a_m   \right]u^m\\
   &=1+\sum_{m\geq 1}\left( m+1 \right)^2a_mu^m\\  
    &=\left(1+u+\cdots+u^n\right)\left(1+\sum_{m\geq 1} a_mu^m\right)\\
   &=1+\sum_{m\geq 1}\left(\sum_{m\geq l\geq \mathrm{max}\{m-n,\,0\}}a_l  \right)u^m.
   \end{split}
   \]  
   It follows that
   \[
   (m+1)^2a_m=\sum_{l=\mathrm{max}\{m-n,\,0\}}^m a_l, \;m\geq 1.
   \]
   By combining similar terms, one has 
   \[
   a_0=1, a_{m+1}=\frac{1}{(m+1)(m+3)}\left( \sum_{l=\mathrm{max}\{m+1-n,\,0\}}^m a_l  \right),  m\geq 0. \tag{2}
   \]
   Define  \[
   b_0=1, b_{m+1}=\frac{1}{(m+1)(m+3)}\left( \sum_{l=0}^m b_l  \right),  m\geq 0.
   \]
      One can check that 
   \[
   b_m=\frac{2}{(m+1)(m+2)},  \;m\geq 0.   \]
   Thus by mathematic induction, one has
   \[
   a_m=b_m=\frac{2}{(m+1)(m+2)}, \;\;0\leq m\leq n, \tag{3}
   \]
   and 
   \[
   0<a_m<b_m=\frac{2}{(m+1)(m+2)},\;\; m>n.    \tag{4}
   \]
   By the formulas $(2), (3)$ and $(4)$, for $m\geq n$, we have 
   \[
   \begin{split}
  0&<a_{m+1}\\
  &\leq \frac{1}{(m+1)(m+3)} \sum_{l=m+1-n}^m \frac{2}{(l+1)(l+2)} \\
  &\leq  \frac{2}{(m+1)(m+3)}  \left( \frac{1}{m+2-n}-\frac{1}{m+2}    \right)\\
    \end{split}
    \]
    As a result,  \[
    0<a_m\leq \frac{2}{m(m+2)}\left(  \frac{1}{m+1-n}-\frac{1}{m+1} \right),\;\;\;m>n.
    \]

     $\hfill\Box$\\   
     
     \begin{lem}\label{fra}
     Define 
     \[
     E_2=\Big{\{}   \frac{1}{2^{k_1}}+\frac{1}{2^{k_1+k_2}}+\cdots+\frac{1}{2^{k_1+\cdots+k_r}}+\cdots   \Big{|}\;k_i\geq 2, i\geq 1  \Big{\}}.
     \]
     Then $    \mathrm{dim}_{\mathrm{H}}\,E_2=  \frac{\mathrm{log}\; \frac{1+\sqrt{5}}{2}}{\mathrm{log}\;2}  .   $
     \end{lem}
     \noindent{\bf Proof:} For $k\geq 2$, define $$\alpha_k: [0,1]\rightarrow [0,1],$$
     \[
     x\mapsto \frac{1}{2^k}+\frac{x}{2^k}.
     \]
     It is easy to check that 
     \[
     E_2=\bigcup_{k\geq 2} \alpha_k(E_2), \;\; \alpha_{k_1}(E_2)  \bigcap  \alpha_{k_2}(E_2)=\emptyset, \;k_1\neq k_2 \]
     and the $\alpha_k,k\geq 2$ satisfy the open set condition by taking $V$ as the open interval $(0,1)$. By theorem $9.3$ in \cite{fal}, the Hausdorff dimention of $E_2$ is given by 
     \[
     \sum_{k\geq 2}\frac{1}{2^{ks}}=\frac{1}{2^s(2^s-1)}=1.
     \]
     So we have 
     \[
     \mathrm{dim}_{\mathrm{H}}\,E_2=s=  \frac{\mathrm{log}\; \frac{1+\sqrt{5}}{2}}{\mathrm{log}\;2}  .      \]
               $\hfill\Box$\\

     Now we are ready to prove partial results of Theorem \ref{rat}.
     
    $(i)$ By definition, we have $\eta_n(\widehat{\mathcal{T}})\subseteq (1, +\infty)$ and 
     \[
     \eta_n((1,1,\cdots,1\cdots))=f_n(1)=\mathrm{e}^{1+\frac{1}{2}+\cdots+\frac{1}{n}} .    \]
     For $(k_1,k_2,\cdots, k_r,\cdots)\in \widehat{\mathcal{T}}$, one has 
     \[
     \begin{split}
    &\;\;\;\; \mathop{\int}_{[0,1]^{k_1+\cdots+k_{r-1}}}  \left[1+t_n(x_1\cdots x_{k_1})\right]\left[1+t_n(x_1\cdots x_{k_1+k_2})\right]\cdots \left[ 1+t_n(x_1\cdots x_{k_1+\cdots+k_{r-1}})\right]\\
    & \cdot dx_1dx_2\cdots dx_{k_1+\cdots+k_{r-1}}<  \mathop{\int}_{[0,1]^{k_1+\cdots+k_{r}}}   \left[1+t_n(x_1\cdots x_{k_1})\right]\left[1+t_n(x_1\cdots x_{k_1+k_2})\right]\cdots\\
    &\cdot  \left[ 1+t_n(x_1\cdots x_{k_1+\cdots+k_{r}})\right] dx_1dx_2\cdots dx_{k_1+\cdots+k_{r}}   \\  \end{split}  \]
     and 
     \[
     \begin{split}
     & \mathop{\int}_{[0,1]^{k_1+\cdots+k_{r}}}   \left[1+t_n(x_1\cdots x_{k_1})\right]\left[1+t_n(x_1\cdots x_{k_1+k_2})\right]\cdots  \left[ 1+t_n(x_1\cdots x_{k_1+\cdots+k_{r}})\right] \\
    & \cdot dx_1dx_2\cdots dx_{k_1+\cdots+k_{r}}  \leq   \mathop{\int}_{[0,1]^{k_1+\cdots+k_{r}}}   \left[1+t_n(x_1)\right]\left[1+t_n(x_1 x_{2})\right]\cdots  \left[ 1+t_n(x_1\cdots x_{k_1+\cdots+k_{r}})\right]\\  
    &\cdot dx_1dx_2\cdots dx_{k_1+\cdots+k_{r}}  .\end{split}
    \]
    In conclusion, the limit $\eta_n\left( (k_1,\cdots,k_r,\cdots)  \right)$ exists and    $$\eta_n(\widehat{\mathcal{T}})\subseteq (1, \mathrm{e}^{1+\frac{1}{2}+\cdots+\frac{1}{n}} ].$$
    The statement $\eta_n(\widehat{\mathcal{T}}-\widehat{\mathcal{T}}^R)= (1, \mathrm{e}^{1+\frac{1}{2}+\cdots+\frac{1}{n}} ]$  will be proved in the next section.
    
    $(ii)$ Firstly we will prove that: for $u\in (0,1)$ and $n\geq 1$,
    \[
     \int^1_0 \left[1+t_n(ux)\right]g_n(ux)dx <1+u+\cdots+u^n  .
    \]
     By Lemma \ref{two}, one has 
     \[
     \begin{split}
     &\;\;\;\;  \int^1_0 \left[1+t_n(ux)\right]g_n(ux)dx     \\
     &\leq   1+\sum_{m= 1}^n \frac{2}{m+2} u^m  + \sum_{m>n}(m+1) \frac{2}{m(m+2)}\left(  \frac{1}{m+1-n}-\frac{1}{m+1} \right) u^m   \\
     &\leq  1+\sum_{m= 1}^n \frac{2}{m+2} u^m  +\left( \sum_{m>n} \frac{2n}{(m+1-n)m(m+2)}\right) u^{n+1} \end{split}
     \]
     For $n=1$, 
       \[
     \begin{split}
     &\;\;\;\;  \int^1_0 \left[1+t_1(ux)\right]g_1(ux)dx     \\
         &\leq  1+\frac{2}{3} u  +\left( \sum_{m>1} \frac{2}{m^2(m+2)}\right) u^{2} \\
         &<1+\left( \sum_{m\geq 1} \frac{2}{m^2(m+2)}\right) u .        \end{split}
     \]
     It is easy to check that 
     \[
     \sum_{m\geq 1} \frac{2}{m^2(m+2)}= \zeta(2)-\frac{3}{4}=0.894934\cdots<1.   \]
     So
      \[
          \int^1_0 \left[1+t_1(ux)\right]g_1(ux)dx     <1+u.
          \]     
     For $n\geq 2$, one has 
      \[
     \begin{split}
     &\;\;\;\;  \int^1_0 \left[1+t_n(ux)\right]g_n(ux)dx     \\
     &\leq  1+\sum_{m= 1}^n \frac{2}{m+2} u^m  +\left( \sum_{m>n} \frac{2n}{(m+1-n)m(m+2)}\right) u^{n+1}\\
     \end{split}
     \]
      and 
      \[
      \begin{split}
      &\;\;\;\; \sum_{m>n} \frac{2n}{(m+1-n)m(m+2)}\\
      &=   \sum_{m>n}\left[   \frac{1}{n+1}\left(\frac{1}{m+1-n}-\frac{1}{m+2}\right)+ \frac{1}{n-1}\left(\frac{1}{m+1-n}-\frac{1}{m}\right)  -\left(\frac{1}{m}-\frac{1}{m+2}     \right)   \right]         \\
      &=\frac{1}{n+1}\left( \frac{1}{2}+\frac{1}{3}+\cdots +\frac{1}{n+2}      \right)+\frac{1}{n-1}\left( \frac{1}{2}+\cdots +\frac{1}{n}      \right)   -\left(\frac{1}{n+1}+\frac{1}{n+2}\right)      \\
      &<\frac{1}{2}+\frac{1}{2}- \left(\frac{1}{n+1}+\frac{1}{n+2}\right) . \end{split}
      \]
      Thus 
       \[
     \begin{split}
     &\;\;\;\;  \int^1_0 \left[1+t_n(ux)\right]g_n(ux)dx     \\
     &\leq  1+\sum_{m= 1}^n \frac{2}{m+2} u^m  +\left[ 1 - \left(\frac{1}{n+1}+\frac{1}{n+2}\right)  \right] u^{n+1}\\
     &\leq 1+\sum_{m= 1}^{n-1} \frac{2}{m+2} u^m  +\frac{n}{n+1}u^n\\
     &<1+u+\cdots+u^n.
        \end{split}
     \]
           
    For $(k_1,k_2,\cdots,k_r,\cdots)\in  \widehat{\mathcal{T}}$, by lemma \ref{two} we have 
    \[
    \begin{split}
    &\;\;\;\;\eta_n\left( (k_1,\cdots, k_{r-1},k_{r}+1, \{2\}^{\infty})   \right)\\
    &= \mathop{\lim}_{l\rightarrow +\infty}  T_n(k_1,\cdots,k_{r-1}, k_r+1,\{2\}^{l})    \\
    &=\int_{[0,1]^{k_1+\cdots+k_r+1}} \left[ 1+t_n(x_1\cdots x_{k_1})  \right]\cdots \left[  1+t_n(x_1\cdots x_{k_1+\cdots+k_{r-1}} ) \right] \left[  1+t_n(x_1\cdots x_{k_1+\cdots+k_{r}+1} ) \right]    \\
    &\;\;\;\;\cdot g_n(x_1\cdots x_{k_1+\cdots +k_r+1})  dx_1\cdots dx_{k_1+\cdots +k_r+1}\\
    &<\int_{[0,1]^{k_1+\cdots+k_r}} \left[ 1+t_n(x_1\cdots x_{k_1})  \right]\cdots \left[  1+t_n(x_1\cdots x_{k_1+\cdots+k_{r-1}} ) \right]\left[1+t_n\left(x_1\cdots x_{k_1+\cdots +k_r}\right)\right] \\
    &\;\;\;\;\cdot dx_1\cdots dx_{k_1+\cdots +k_r}\\
    &\leq T_n(k_1,\cdots,k_{r-1}, k_r)  \\
    &< \eta_n\left((k_1,k_2,\cdots,k_r,\cdots)\right). \end{split}
    \]   
     For ${\bf k,l}\in  \widehat{\mathcal{T}}_2$ and     $${\bf k}=(k_1,\cdots,k_r\cdots)\succ {\bf l}=(l_1,\cdots,l_r,\cdots),$$there is an $r_0\geq 1$ such that 
     \[
     k_i=l_i, i\leq r_0-1,\; k_{r_0}<l_{r_0}.
     \]
     We have 
     \[
     \eta_n({\bf l})\leq \eta_n((l_1,\cdots,l_{r_0-1},l_{r_0},\{2\}^{\infty}))<T_n(l_1,\cdots, l_{r_0-1}, k_{r_0})<\eta_{n}({\bf k}).
     \]
     Similarly, if ${\bf l}\succ {\bf k}$, then $\eta_n({\bf l})> \eta_n({\bf k})$. In conclusion, for ${\bf k,l}\in  \widehat{\mathcal{T}}_2$, we have $\eta_n({\bf k})> \eta_n({\bf l})$ if and only if    ${\bf k}\succ {\bf l}$. 
     
     $(iii)$  For ${\bf k,l}\in  \widehat{\mathcal{T}}_2$ and     $${\bf k}=(k_1,\cdots,k_r\cdots)\succ {\bf l}=(l_1,\cdots,l_r,\cdots),$$there is an $r_0\geq 1$ such that 
     \[
     k_i=l_i, i\leq r_0-1,\; k_{r_0}<l_{r_0}.
     \] 
    We want to show that for some fixed $c_1$ and $c_2$, 
     \[
    \frac{c_1}{2^{k_1+\cdots +k_{r_0-1}+k_{r_0}}}  <\eta_{n}({\bf k})-  \eta_n({\bf l})<  \frac{c_2}{2^{k_1+\cdots +k_{r_0-1}+k_{r_0}}}.\]
    We have 
    \[
    \begin{split}
    &\;\;\;\;\eta_{n}({\bf k})-  \eta_n({\bf l})\\
    &< \eta_n((k_1,\cdots,k_{r_0-1},k_{r_0},\{2\}^{\infty}))   - T_n(l_1,\cdots,l_{r_0-1},l_{r_0})    \\
    &< \eta_n((k_1,\cdots,k_{r_0-1},k_{r_0},\{2\}^{\infty}))   - T_n(l_1,\cdots,l_{r_0-1})  \\
    &\leq  \int_{[0,1]^{k_1+\cdots+k_{r_0-1}}} \left[  1+t_n(x_1\cdots x_{k_1})  \right] \left[ 1+t_n(x_1\cdots x_{k_1+\cdots +k_{r_0-1}})\right]\\
     &\cdots\left[ ( 1+t_n(x_1\cdots x_{k_1+\cdots +k_{r_0}}))  g_n(x_1\cdots x_{k_1+\cdots+k_{r_0}})-1\right] dx_1\cdots dx_{k_1+\cdots+k_{r_0}} \\
    &\leq  \int_{[0,1]^{k_1+\cdots+k_{r_0}}} \left[  1+t_n(x_1\cdots x_{k_1})  \right]  \cdots \left[ 1+t_n(x_1\cdots x_{k_1+\cdots +k_{r_0-1}})\right] \\
    &\cdot \sum_{m\geq 1}    \frac{ (m+1)a_m}{(m+1)^{k_{r_0}}}   \left(  x_1\cdots x_{k_1+\cdots +k_{r_0-1}}   \right)^m     dx_1\cdots dx_{k_1+\cdots+k_{r_0-1}} \\
   &< \frac{n}{ 2^{k_{r_0}} }\int_{[0,1]^{k_1+\cdots+k_{r_0-1}}} \left[  1+nx_1\cdots x_{k_1}  \right]  \cdots \left[ 1+nx_1\cdots x_{k_1+\cdots +k_{r_0-1}}\right] x_1\cdots x_{k_1+\cdots+k_{r_0-1}}\\
    &\cdot dx_1\cdots dx_{k_1+\cdots+k_{r_0-1}} \\
    &<\frac{n e^{2n}}{2^{k_1+\cdots+k_{r_0}}}.
     \end{split}
    \]
    Here the last inequality follows from Lemma \ref{leq}.
    
    On the other hand,
      \[
    \begin{split}
    &\;\;\;\;\eta_{n}({\bf k})-  \eta_n({\bf l})\\
    &>T_n(l_1,\cdots, l_{r_0-1}, k_{r_0}) -\eta_n((l_1,\cdots,l_{r_0-1},k_{r_0}+1,\{2\}^{\infty}))\\
 &\geq  \int_{[0,1]^{k_1+\cdots+k_{r_0}}} \left[  1+t_n(x_1\cdots x_{l_1})  \right]  \cdots \left[ 1+t_n(x_1\cdots x_{l_1+\cdots +l_{r_0-1}})\right]\\
 &\cdot  \Big{[}1+t_n(x_1\cdots x_{l_1+\cdots+l_{r_0-1}+k_{r_0}})-\int^1_0 \left[ 1+t_n( ux_1\cdots x_{l_1+\cdots+l_{r_0-1}+k_{r_0}})  \right] \\
 &\;\;\;\;g_n(u x_1\cdots x_{l_1+\cdots+l_{r_0-1}+k_{r_0}})du\Big{]} dx_1\cdots dx_{k_1+\cdots+k_{r_0}} \\
&\geq c_1 \int_{[0,1]^{k_1+\cdots+k_{r_0}}}x_1\cdots x_{l_1+\cdots+l_{r_0-1}+k_{r_0}}dx_1\cdots dx_{k_1+\cdots+k_{r_0}} \\
&\geq \frac{c_1}{2^{k_1+\cdots+k_{r_0-1}+k_{r_0}}}.\end{split} 
    \]        
    Here 
    \[
    c_1= \mathrm{min}\bigg{\{} \frac{7}{4}-\zeta(2)  , \frac{1}{3}\bigg{\}}.
    \]
    From the above analysis, we have 
      \[
    \frac{c_1}{2^{k_1+\cdots +k_{r_0-1}+k_{r_0}}}  <\eta_{n}({\bf k})-  \eta_n({\bf l})<  \frac{c_2}{2^{k_1+\cdots +k_{r_0-1}+k_{r_0}}}\tag{2}
    \]   
     for $c_1$ and $c_2=ne^{2n}$.
     
     Define $$i_n: E_2\rightarrow \eta_n(\widehat{\mathcal{T}}_2),$$
     \[
     \frac{1}{2^{k_1}}+\frac{1}{2^{k_1+k_2}}+\cdots+\frac{1}{2^{k_1+\cdots+k_r}}+\cdots \mapsto \eta_n((k_1,k_2,\cdots,k_r,\cdots)).
     \] From the formula $(2)$  it is easy to check that $i_n$ is bi-Lipschitz equivalent. By Lemma \ref{fra}, we have 
     \[
     \mathrm{dim}_{\mathrm{H}}\,\eta_n(\widehat{\mathcal{T}}_2)=  \mathrm{dim}_{\mathrm{H}}\,E_2=  \frac{\mathrm{log}\; \frac{1+\sqrt{5}}{2}}{\mathrm{log}\;2}  .  \]
     
     \begin{rem}\label{sleq}
     In fact, from the simple observation that for $k_1,\cdots,k_r)\neq (1,\cdots,1)$,
     \[
     (1+x_1\cdots x_{k_1}) (1+x_1x_2\cdots x_{k_1})  \cdots (1+x_1\cdots x_{k_1+\cdots+k_r})  <   \prod_{i=1}^{k_1+\cdots+k_r}  (1+x_1\cdots x_i), \]
     one can deduce that 
     \[\eta_n\left((k_1,k_2,\cdots,k_r, \cdots) \right)<e^{1+\frac{1}{2}+\cdots+\frac{1}{n}}\] 
     if $(k_1,k_2,\cdots,k_r, \cdots)\neq (\{1\}^\infty) $ .   \end{rem}
     
     \section{Function decomposition}
     In this section we will study the function $F_n=\eta_n\circ \tau^{-1}$. Although $F_n$ is not continuous at every point, we will show that it is close to a continuous function in some sense.
     
     \begin{lem}\label{exp} Define $\Delta_n(u)$ as 
     \[
    \Delta_n(u)= \mathrm{exp}\left(u+\frac{u^2}{2}+\cdots+\frac{u^n}{n}    \right)-\left( 1+u+\cdots+u^n     \right) .    \]
     For $u\in (0,1)$, 
     \[
    \frac{n\cdot  u^{n+1}}{n+1} <\Delta_n(u)<2n\cdot u^{n+1}.
     \]
     \end{lem}
     \noindent{\bf Proof:}
     Define 
     \[
     \Delta_n(u)=\sum_{m\geq 1} c_m u^m.
     \]
     It is clear that $c_m>0$ for $m>n$.
     Since 
     \[
     \begin{split}
     &\;\;\;\;\mathop{\lim}_ {u\rightarrow 0^+} \frac{\Delta_n(u)}{u^{n+1}} \\
     &=\mathop{\lim}_ {u\rightarrow 0^+}  \frac{\mathrm{exp}\left(\mathrm{log}\,\frac{1}{1-u}- \sum\limits_{m\geq n+1} \frac{u^m}{m}   \right)-\left( 1+u+\cdots+u^n     \right)  }{u^{n+1}}   \\
     &= \mathop{\lim}_ {u\rightarrow 0^+}  \frac{\mathrm{exp}\left(  -\sum\limits_{m\geq n+1}\frac{u^m}{m} \right)-(1-u^{n+1}) }{(1-u)u^{n+1}}\\
     &=\frac{n}{n+1},\end{split}
     \]
     one has $c_m=0$, $m\leq n$ and $c_{n+1}=\frac{n}{n+1}$. As a result, $\Delta_n(u)>\frac{n}{n+1}\cdot  u^{n+1} $.
     As $\frac{\Delta_n(u)}{u^{n+1}}=\sum\limits_{m\geq n+1}c_m u^{m-n-1} $ is  a  strictly increasing function, one has
     \[
     \begin{split}
     &\;\;\;\;  \frac{\Delta_n(u)}{u^{n+1}}      \\
     &<\Delta_n(1)            \\
     &\leq  \mathrm{exp}\left(1+\frac{1}{2}+\cdots+\frac{1}{n}    \right)-(n+1) \\
     &<  \mathrm{exp}\left(1+\int^2_1\frac{du}{u}+\cdots+\int^n_{n-1}\frac{du}{u}    \right)-(n+1)  \\
     &\leq en-n-1\\
     &<2n. \end{split}
     \]   
     So we have 
     \[
     \Delta_n(u)<2n\cdot u^{n+1}.     \]      $\hfill\Box$\\         
     
        \begin{lem}\label{del}
       For $k_1,\cdots,k_r\geq 1$ and 
     \[
     y=\tau\left((k_1,\cdots,k_{r-1},k_r+1,\{1\}^\infty)\right)=\frac{1}{2^{k_1}}+\frac{1}{2^{k_1+k_2}}+\cdots+\frac{1}{2^{k_1+\cdots+k_r}},     \]  
     define 
     \[
     \begin{split}
    & h_n(y)=\mathop{\int}_{[0,1]^{k_1+\cdots+k_r}}\left[ 1+t_n(x_1\cdots x_{k_1})\right] \cdots \left[   1+t_n(x_1\cdots x_{k_1+\cdots +k_{r-1}})  \right]\Delta_n(x_1\cdots x_{k_1+\cdots+k_r})\\
    &\;\;\;\;\;\;\;\;\;\;\;\;\cdot dx_1\cdots dx_{k_1+\cdots+k_r}.\\
     \end{split}
     \]
     Then 
     \[
             \frac{\frac{n}{n+1}}{(n+2)^{k_1+\cdots+k_r}}\leq     h_n(y)\leq \frac{2n e^{(n+2)n}}{(n+2)^{k_1+\cdots+k_r}}.
     \]
     \end{lem}
      \noindent{\bf Proof:}
      By Lemma \ref{exp} and Lemma \ref{leq}, one has 
      \[
      \begin{split}
      &\;\;\;\;  h_n(y)    \\
      &\leq 2n   \mathop{\int}_{[0,1]^{k_1+\cdots+k_r}}\left[ 1+t_n(x_1\cdots x_{k_1})\right] \cdots\left[   1+t_n(x_1\cdots x_{k_1+\cdots +k_{r-1}})  \right](x_1\cdots x_{k_1+\cdots+k_r})^{n+1}\\
      &\;\;\;\;\;\;\cdot  dx_1\cdots dx_{k_1+\cdots+k_r}    \\
      &< 2n   \mathop{\int}_{[0,1]^{k_1+\cdots+k_r}}\left[ 1+nx_1\cdots x_{k_1}\right] \cdots \left[   1+nx_1\cdots x_{k_1+\cdots +k_{r-1}}  \right](x_1\cdots x_{k_1+\cdots+k_r})^{n+1}\\
      &\;\;\;\;\;\;\cdot  dx_1\cdots dx_{k_1+\cdots+k_r}    \\   
      &\leq \frac{2n e^{(n+2)n}}{(n+2)^{k_1+\cdots+k_r}}  .               \\   \end{split}
      \]
      By Lemma \ref{exp}, one has      
      \[
      \begin{split}
      &\;\;\;\;\; h_n(y)\\
      &> \mathop{\int}_{[0,1]^{k_1+\cdots+k_r}}\Delta_n(x_1\cdots x_{k_1+\cdots+k_r})
    dx_1\cdots dx_{k_1+\cdots+k_r}   \\
    &>\frac{n}{n+1} \mathop{\int}_{[0,1]^{k_1+\cdots+k_r}}(x_1\cdots x_{k_1+\cdots+k_r})^{n+1}
    dx_1\cdots dx_{k_1+\cdots+k_r}\\
    &\geq   \frac{n}{n+1} \frac{1}{(n+2)^{k_1+\cdots+k_r}}  .   \end{split}
      \]
     $\hfill\Box$\\     
     
     Define $$H_n: (0, 1]\rightarrow (0, +\infty),$$    
     \[
     x\mapsto \sum_{\substack{y<x\\ y\in \tau \left(\widehat{\mathcal{T}}^R \right)}} h_n(y).
     \] 
     In order to show that the map $H_n$ is well-defined, it suffices to prove that 
     \begin{lem}\label{wel}
     The series $$\sum_{ y\in \tau \left(\widehat{\mathcal{T}}^R \right)} h_n(y)$$ is absolutely convergent.
          \end{lem}
  \noindent{\bf Proof:} By Lemma \ref{del}, one has 
  \[
  \begin{split}
  &\;\;\;\; \sum_{ y\in \tau \left(\widehat{\mathcal{T}}^R \right)} h_n(y) \\
  &= \sum_{r\geq 1} \sum_{k_1,\cdots,k_r\geq 1}  \mathop{\int}_{[0,1]^{k_1+\cdots+k_r}}\left[ 1+t_n(x_1\cdots x_{k_1})\right] \cdots \left[   1+t_n(x_1\cdots x_{k_1+\cdots +k_{r-1}})  \right]\\
    &\;\;\;\;\;\;\;\;\;\;\;\;\cdot \Delta_n(x_1\cdots x_{k_1+\cdots+k_r})  dx_1\cdots dx_{k_1+\cdots+k_r}  \\
    &\leq 2n e^{(n+2)n}  \sum_{r\geq 1} \sum_{k_1,\cdots,k_r\geq 1}    \frac{1}{(n+2)^{k_1+\cdots+k_r}}\\
    &\leq 2n e^{(n+2)n}  \sum_{r\geq 1} \left[ \sum_{k\geq 1}    \frac{1}{(n+2)^{k}}\right]^r \\
    & \leq 2 e^{(n+2)n} . \end{split}
  \]
  $\hfill\Box$\\ 

Since the set $\tau(\widehat{\mathcal{T}}^{R})$ is dense in $(0, 1]$, the function $H_n$ is strictly increasing. One can also check that $H_n$ is continuous at $(0,1]-\tau(\widehat{\mathcal{T}}^{R})$ and discontinuous at $\tau(\widehat{\mathcal{T}}^{R})$.   
\begin{lem}\label{hau}
For the structure of the image  $\mathrm{Im}\,(H_n)$,  one has that 
$\mathrm{Im}\,(H_n)$ is bi-Lipschitz      equivalent to    $C_{\frac{1}{n+2}}-  C_{\frac{1}{n+2}}^b$    and   \[  \mathrm{dim}_{\mathrm{H}} \mathrm{Im}\,(H_n)=\frac{\mathrm{log}\,2}{\mathrm{log}(n+2)} .\]\end{lem}  
   \noindent{\bf Proof:}
   By definition, one has 
   \[
   C_{\frac{1}{n+2}}-  C_{\frac{1}{n+2}}^b=\Big{\{}\frac{n+1}{(n+2)^{k_1}}+\frac{n+1}{(n+2)^{k_1+k_2}} +\cdots+ \frac{n+1}{(n+2)^{k_1+\cdots+k_r}} +\cdots \,\Big{|}\, k_i\geq 1, i\geq 1\Big{\}}.\]
   Define \[j_n:   C_{\frac{1}{n+2}}- C_{\frac{1}{n+2}}^b  \rightarrow    \mathrm{Im}\,(H_n),\] 
      \[
   \frac{n+1}{(n+2)^{k_1}}+\frac{n+1}{(n+2)^{k_1+k_2}} +\cdots+ \frac{n+1}{(n+2)^{k_1+\cdots+k_r}} +\cdots  \mapsto H_n(\tau(k_1,k_2,\cdots, k_r,\cdots)).  \]
   It is easy to check that $j_n$ is bijective.
   For ${\bf k,l}\in  \widehat{\mathcal{T}}$ and     $${\bf k}=(k_1,\cdots,k_r\cdots)\succ {\bf l}=(l_1,\cdots,l_r,\cdots),$$there is an $r_0\geq 1$ such that 
     \[
     k_i=l_i, i\leq r_0-1,\; k_{r_0}<l_{r_0}.
     \] 
     We have 
     \[
     H_n\left( \tau({\bf k})\right)-H_n\left( \tau({\bf l})\right)  
     =\begin{cases}
     \sum_{\substack{   {\bf p} \in \widehat{\mathcal{T}}^R      \\    \bf k \succ \bf p  \succ \bf l  }}  h_n\left(\tau({\bf p})\right)+  h_n\left(\tau({\bf l})\right),   & {\bf l}\in\widehat{\mathcal{T}}^R;\\
   \sum_{\substack{   {\bf p} \in \widehat{\mathcal{T}}^R      \\    \bf k \succ \bf p  \succ \bf l  }}  h_n\left(\tau({\bf p})\right), & {\bf l}\notin\widehat{\mathcal{T}}^R .\\     \end{cases}  
         \]
         It follows that 
        \[
        \begin{split}
        &\;\;\;\;  H_n\left( \tau({\bf k})\right)-H_n\left( \tau({\bf l})\right)    \\
               &> \sum_{l_{r_0}-1\geq p_{r_0}\geq k_{r_0}} h_n\left( \frac{1}{2^{k_1}}+\cdots+\frac{1}{2^{k_1+\cdots+k_{r_0-1}}}+\frac{1}{2^{k_1+\cdots+k_{r_0-1}+p_{r_0}}} \right)\\
        &\geq   \frac{\frac{n}{n+1}}{(n+2)^{k_1+\cdots+k_{r_0}}} .       \end{split}
        \]
        Since $l_{r_0}>k_{r_0}\geq 1$, by Lemma \ref{del}, if ${\bf l}\in \widehat{\mathcal{T}}^R$,  one has
        \[
        h_n\left( \tau({\bf l})    \right) \leq \frac{2n e^{(n+2)n}}{(n+2)^{k_1+\cdots+k_{r_0}}}        \]
        
         As a result, one has 
         \[
         \begin{split}
         & \;\;\;\;  H_n\left( \tau({\bf k})\right)-H_n\left( \tau({\bf l})\right)    \\
        &\leq     \sum_{\substack{   {\bf p} \in \widehat{\mathcal{T}}^R      \\    \bf k \succ \bf p  \succ \bf l  }}  h_n\left(\tau({\bf p})\right) +\frac{2n e^{(n+2)n}}{(n+2)^{k_1+\cdots+k_{r_0}}}        \\
        &\leq  \sum_{k_{r_0}\leq p_{r_0}\leq l_{r_0}} \Bigg{[}h_n\Bigg{(}\frac{1}{2^{k_1}}+\cdots+\frac{1}{2^{k_1+\cdots+k_{r_0-1}}}+ \frac{1}{2^{k_1+\cdots+k_{r_0-1}+p_{r_0}}}      \Bigg{)}
        + \sum_{m\geq 1}\sum_{q_1,\cdots,q_m\geq 1} \\
        &h_n\Bigg{(}  \frac{1}{2^{k_1}}+\cdots+\frac{1}{2^{k_1+\cdots+k_{r_0-1}}}+ \frac{1}{2^{k_1+\cdots+k_{r_0-1}+p_{r_0}}}   +\frac{1}{2^{k_1+\cdots+k_{r_0-1}+p_{r_0}+q_1}} \\
        &+\cdots+ \frac{1}{2^{k_1+\cdots+k_{r_0-1}+p_{r_0}+q_1+\cdots+q_m}}    \Bigg{) }\Bigg{]} +\frac{2n e^{(n+2)n}}{(n+2)^{k_1+\cdots+k_{r_0}}}   \\
        &\leq 2n e^{(n+2)n}\Bigg{[} \frac{1}{(n+2)^{k_1+\cdots+k_{r_0}}} +\sum_{k_{r_0}\leq p_{r_0}\leq l_{r_0}} \Bigg{(} \frac{1}{(n+2)^{k_1+\cdots+k_{r_0-1}+p_{r_0}}}+\sum_{m\geq 1}\sum_{q_1,\cdots,q_m\geq 1} \\
        &   \frac{1}{(n+2)^{k_1+\cdots+k_{r_0-1}+p_{r_0}+q_1+\cdots+q_m } } \Bigg{)} \Bigg{]}  \\
        &\leq 2n e^{(n+2)n}\Bigg{[} \frac{1}{(n+2)^{k_1+\cdots+k_{r_0}}} +\sum_{k_{r_0}\leq p_{r_0}\leq l_{r_0}} \Bigg{(} \frac{1}{(n+2)^{k_1+\cdots+k_{r_0-1}+p_{r_0}}}+\sum_{m\geq 1} \\
        &   \frac{1}{(n+2)^{k_1+\cdots+k_{r_0-1}+p_{r_0} } }\cdot  \frac{1}{(n+1)^m} \Bigg{)} \Bigg{]}  \\   
           \end{split}
        \]
        \[
        \begin{split}
        &  \leq 2n e^{(n+2)n}\Bigg{[} \frac{1}{(n+2)^{k_1+\cdots+k_{r_0}}} +\sum_{k_{r_0}\leq p_{r_0}\leq l_{r_0}}  \frac{1}{(n+2)^{k_1+\cdots+k_{r_0-1}+p_{r_0}}}\cdot\left(1+\frac{1}{n}\right)\Bigg{]} \\ 
               &  \leq 2n e^{(n+2)n}\Bigg{[} \frac{1}{(n+2)^{k_1+\cdots+k_{r_0}}} +  \frac{2}{(n+2)^{k_1+\cdots+k_{r_0-1}+k_{r_0}}}\cdot\left(1+\frac{1}{n}\right)\Bigg{]}       \\ 
        &\leq  \frac{(6n+4)e^{(n+2)n} }{(n+2)^{k_1+\cdots+k_{r_0}}}.  \\ 
        \end{split}
        \]
        
        In conclusion, the map $j_n$ is bi-Lipschitz equivalent. Thus one has 
        \[
         \mathrm{dim}_{\mathrm{H}} \mathrm{Im}\,(H_n)=\mathrm{dim}_{\mathrm{H}}   \left( C_{\frac{1}{n+2}}- C_{\frac{1}{n+2}}^b\right)  = \mathrm{dim}_{\mathrm{H}}   \left( C_{\frac{1}{n+2}}\right)=\frac{\mathrm{log}\,2}{\mathrm{log}(n+2)} .        \]        
     $\hfill\Box$\\ 
     
         From the above analysis in this section, Theorem \ref{dec} is reduced to the following theorem.
         \begin{Thm}\label{equ}
         Define $G_n$ as 
         \[
         G_n: (0,1]\rightarrow (1,+\infty),
         \]
         \[
         x\mapsto F_n(x)+H_n(x).
         \]
         Then $G_n$ is  continuous on $(0, 1]$ and strictly increasing. 
         \end{Thm}
          \noindent{\bf Proof:}   For $x_0\in (0,1)-\tau(\widehat{\mathcal{T}}^{R})$, we  want to show that $F_n$ is continuous at $x_0$. We have 
          \[
          x_0=\frac{1}{2^{k_1}}+\frac{1}{2^{k_1+k_2}}+\cdots+\frac{1}{2^{k_1+\cdots+k_r}}+\cdots
          \]
           and $k_r\geq 2 $ for infinitely many $r$. For  a sequence $\{y_m\in (0,1)\,|\, m\geq 1\}$ which satisfies 
           \[
           x_0=\mathop{\lim}_{m\rightarrow +\infty} y_m,
           \]
           \[
           y_m =  \tau\left((l_{m,1},l_{m,2},\cdots, l_{m,r},\cdots  )     \right),          \]
              there is a sequence $\{a_r\in \mathbb{N}\,|\, r\geq 1\}$ such that 
              \[
              \mathop{\lim}_{n\rightarrow +\infty} a_r=+\infty,              \]
              and 
              \[
              (l_{m,1}, l_{m,2},\cdots, l_{m,r})=(k_1,k_2,\cdots,k_r), \;\;m\geq a_r, r\geq 1.
              \] 
              As a result, for $m\geq a_r$
              \[
             T_n(   k_1,k_2,\cdots,k_{r}       ) <F_n(y_m)< F_n\left(\tau((k_1,k_2,\cdots,k_{r}, \{1\}^{\infty}))  \right) .           \]
             For $r\geq 1$, by Lemma \ref{111},  one has 
             \[
             \begin{split}
             &\;\;\;\;F_n\left(\tau((k_1,k_2,\cdots,k_{r}, \{1\}^{\infty}))  \right) -     T_n(   k_1,k_2,\cdots,k_{r}       )    \\
             &=  \mathop{\int}_{[0,1]^{k_1+\cdots+k_r}}   \left[  1+t_n(x_1\cdots x_{k_1})     \right] \cdots \left[    1+t_n(x_1\cdots x_{k_1+\cdots+k_r})      \right] \\
             &\cdot \left[  \mathrm{exp}\left(x_1\cdots x_{k_1+\cdots+k_r}+\frac{ (x_1\cdots x_{k_1+\cdots+k_r} )^2}{2}+\cdots+ \frac{ (x_1\cdots x_{k_1+\cdots+k_r} )^n}{n}    \right)-1   \right] \\
             &\cdot dx_1dx_2\cdots dx_{k_1+\cdots+k_r}     \\
                                &< \mathop{\int}_{[0,1]^{k_1+\cdots+k_r}}   \left[  1+t_n(x_1\cdots x_{k_1})     \right] \cdots \left[    1+t_n(x_1\cdots x_{k_1+\cdots+k_r})      \right] \\
             &\;\;\;\;\cdot \left[  \mathrm{exp}\left(\left(1+\frac{1}{2}+\cdots+\frac{1}{n}\right)x_1\cdots x_{k_1+\cdots+k_r}  \right)-1   \right]dx_1dx_2\cdots dx_{k_1+\cdots+k_r}     \\
               \end{split}
             \]
             \[
             \begin{split}       
              &< \mathop{\int}_{[0,1]^{k_1+\cdots+k_r}}   \left[  1+t_n(x_1\cdots x_{k_1})     \right] \cdots \left[    1+t_n(x_1\cdots x_{k_1+\cdots+k_r})      \right] \\
             &\;\;\;\;\cdot \left[  \mathrm{exp}\left(\left(1+\mathrm{log}\,n\right)x_1\cdots x_{k_1+\cdots+k_r}  \right)-1   \right]dx_1dx_2\cdots dx_{k_1+\cdots+k_r}     \\
              &\leq  \mathop{\int}_{[0,1]^{k_1+\cdots+k_r}}   \left[  1+t_n(x_1\cdots x_{k_1})     \right] \cdots \left[    1+t_n(x_1\cdots x_{k_1+\cdots+k_r})      \right] \\
             &\;\;\;\;\cdot   \sum_{m\geq 1}\frac{(1+\mathrm{log}\, n)^m}{m!} (x_1\cdots x_{k_1+\cdots+k_r})^m   dx_1dx_2\cdots dx_{k_1+\cdots+k_r}     \\      
              &\leq  \sum_{m\geq 1}\frac{(1+\mathrm{log}\, n)^m}{m!} \mathop{\int}_{[0,1]^{k_1+\cdots+k_r}}   \left[  1+t_n(x_1\cdots x_{k_1})     \right] \cdots \left[    1+t_n(x_1\cdots x_{k_1+\cdots+k_r})      \right] \\
             &\;\;\;\;\cdot   x_1\cdots x_{k_1+\cdots+k_r}  dx_1dx_2\cdots dx_{k_1+\cdots+k_r}     \\  
                   &\leq  \sum_{m\geq 1}\frac{(1+\mathrm{log}\, n)^m}{m!} \mathop{\int}_{[0,1]^{k_1+\cdots+k_r}}   \left[  1+t_n(x_1)     \right]\left[  1+t_n(x_1x_2)     \right] \cdots \left[    1+t_n(x_1\cdots x_{k_1+\cdots+k_r})      \right] \\
             &\;\;\;\;\cdot   x_1\cdots x_{k_1+\cdots+k_r}  dx_1dx_2\cdots dx_{k_1+\cdots+k_r}     \\  
              &\leq  \sum_{m\geq 1}\frac{(1+\mathrm{log}\, n)^m}{m!} \mathop{\int}_{[0,1]^{k_1+\cdots+k_r}}   \left[  1+nx_1     \right]\left[  1+nx_1x_2     \right] \cdots \left[    1+nx_1\cdots x_{k_1+\cdots+k_r}    \right] \\
             &\;\;\;\;\cdot   x_1\cdots x_{k_1+\cdots+k_r}  dx_1dx_2\cdots dx_{k_1+\cdots+k_r}     \\      
              &\leq  \sum_{m\geq 1}\frac{(1+\mathrm{log}\, n)^m}{m!} \cdot \frac{e^{2n}}{2^{k_1+k_2+\cdots+k_r}}  .                         \end{split}
                  \]
                  It follows that 
                  \[
                    \mathop{\lim}_{m\rightarrow +\infty} F_n(y_m)=   \mathop{\lim}_{r \rightarrow +\infty} T_n(   k_1,k_2,\cdots,k_{r}       ) =F_n(x_0).              \]
     Thus    $F_n$ is continuous on $(0,1)-\tau(\widehat{\mathcal{T}}^{R})$. By the same analysis one can show that $F_n$ is (left) continuous  at $x_0=1$.

         As the function $H_n$ is continuous on  $(0,1]-\tau(\widehat{\mathcal{T}}^{R})$ and discontinuous at $\tau(\widehat{\mathcal{T}}^{R})$, $G_n$ is continuous on $(0,1]-\tau(\widehat{\mathcal{T}}^{R})$.

         For $$x_0=\frac{1}{2^{k_1}}+\frac{1}{2^{k_1+k_2}}+\cdots+\frac{1}{2^{k_1+\cdots+k_r}}\in \tau(\widehat{\mathcal{T}}^{R}),$$   
         firstly we will show that $G_n$ is left continuous     at $x_0$.    For $\{ y_m\,|\, m\geq 1\}\subseteq (0,x_0)$ and 
         \[
         \mathop{\lim}_{m\rightarrow +\infty} y_m=x_0,
         \]
          \[
           y_m =  \tau\left((p_{m,1},p_{m,2},\cdots, p_{m,r},\cdots  )     \right),          \]          there is a sequence $\{b_l\in \mathbb{N}\,|\, l\geq 1\}$ such that 
              \[
              \mathop{\lim}_{n\rightarrow +\infty} b_l=+\infty,              \]
              and 
              \[(p_{m,1}, p_{m,2},\cdots, p_{m,r}, \cdots, p_{m,r+l})=(k_1,\cdots,k_{r-1},k_r+1,\{1\}^l ), \;\;m\geq b_l, l\geq 1.\]
              Define 
              \[
              z_l=\frac{1}{2^{k_1}}+\cdots+\frac{1}{2^{k_1+\cdots+k_{r-1}} }+\frac{1}{2^{k_1+\cdots+k_{r}+1} }+\cdots+ \frac{1}{2^{k_1+\cdots+k_{r}+l+1} }              \]
              From the above equation, for $m\geq b_l$, one has
              \[
              \begin{split}
              &\;\;\;\;\,  G_n(y_m)\\
              &=F_n(y_m)+H_n(y_m)\\
              &= \eta_n\left((p_{m,1},p_{m,2},\cdots,p_{m,r},\cdots)      \right) +\sum_{\substack{y<y_m\\ y\in \tau \left(\widehat{\mathcal{T}}^R \right)}} h_n(y)   \\
              &\leq \eta_n\left(  (k_1,\cdots,k_{r-1},k_r+1,\{1\}^\infty )     \right) +\sum_{\substack{y<x_0\\ y\in \tau \left(\widehat{\mathcal{T}}^R \right)}} h_n(y)  =G_n(x_0)             \end{split}
              \] and 
               \[
              \begin{split}
              &\;\;\;\;\,  G_n(y_m)\\
              &=F_n(y_m)+H_n(y_m)\\
                      &> T_n\left(  k_1,\cdots,k_{r-1},k_r+1,\{1\}^l   \right) +\sum_{\substack{y< z_l \\ y\in \tau \left(\widehat{\mathcal{T}}^R \right)}} h_n(y)  .             \end{split}
              \] 
              By Lemma \ref{del}, it  follows that 
              \[
              \begin{split}
              &\;\;\;\; G_n(x_0)  -G_n(y_m)         \\
              &\leq  \eta_n\left(  (k_1,\cdots,k_{r-1},k_r+1,\{1\}^\infty )     \right)- T_n\left(  k_1,\cdots,k_{r-1},k_r+1,\{1\}^l   \right) +  \sum_{\substack{z_l\leq y<x_0\\ y\in \tau \left(\widehat{\mathcal{T}}^R \right)}} h_n(y)       \\
              &\leq  \eta_n\left(  (k_1,\cdots,k_{r-1},k_r+1,\{1\}^\infty )     \right)- T_n\left(  k_1,\cdots,k_{r-1},k_r+1,\{1\}^l   \right) \\
              &\;\;\;\;+ 2n e^{(n+2)n}\left(\frac{1}{(n+2)^{k_1+\cdots+k_{r}+l+1} }+ \sum_{m\geq 1} \sum_{l_1,\cdots,l_m\geq 1}  \frac{1}{(n+2)^{k_1+\cdots+k_{r}+l+1+l_1+\cdots+l_m}}   \right)         \\
              &\leq  \eta_n\left(  (k_1,\cdots,k_{r-1},k_r+1,\{1\}^\infty )     \right)- T_n\left(  k_1,\cdots,k_{r-1},k_r+1,\{1\}^l   \right) \\
              &\;\;\;\;+\frac{   2(n+1) e^{(n+2)n} }{(n+2)^{k_1+\cdots+k_{r}+l+1} }.         \\
              \end{split}
              \]
              Since 
              \[
               \eta_n\left(  (k_1,\cdots,k_{r-1},k_r+1,\{1\}^\infty )     \right)=\mathop{\lim}_{l\rightarrow +\infty}T_n\left(  k_1,\cdots,k_{r-1},k_r+1,\{1\}^l   \right)  ,            \]
              we have 
              \[
              G_n(x_0)=\mathop{\lim}_{m\rightarrow+\infty} G_n(y_m).
              \]
              In a word, $G_n$ is left continuous at $x_0$.
              
              Secondly we will show that $G_n$ is right continuous at $x_0$. For $\{ u_m\,|\, m\geq 1\}\subseteq (x_0, 1)$ and 
         \[
         \mathop{\lim}_{m\rightarrow +\infty} u_m=x_0,
         \]
          \[
           u_m =  \tau\left((q_{m,1},q_{m,2},\cdots, q_{m,r},\cdots  )     \right),          \]          there is a sequence $\{c_l\in \mathbb{N}\,|\, l\geq 1\}$ such that 
              \[
              \mathop{\lim}_{n\rightarrow +\infty} c_l=+\infty,              \]
         \[(p_{m,1}, p_{m,2},\cdots, p_{m,r})=(k_1,k_2,\cdots,k_r)           , 
          \]             and  $p_{m,r+1}> l$ for     $m\geq c_l, l\geq 1$.       For  $m\geq c_l, l\geq 1$,       one has  
           \[
           \begin{split}
          &\;\;\;\;G_n(u_m)\\
          &=F_n(u_m)+H_n(u_m)\\
          &\leq F_n\left(\tau((k_1,k_2,\cdots,k_r,l, \{1\}^\infty   )     )\right)+H_n(\tau((k_1,k_2,\cdots,k_r,l, \{1\}^\infty   )     ) ),                        \end{split}
          \]     
                \[
           \begin{split}
          &\;\;\;\;G_n(u_m)\\
          &=F_n(u_m)+H_n(u_m)\\
          &\geq T_n(k_1,k_2, \cdots,k_r)+H_n(x_0 )+h_n(x_0).\\
\end{split}
          \]    
          For convenience, we write $$w_l=  \tau((k_1,k_2,\cdots,k_r,l, \{1\}^\infty   )     ) =\frac{1}{2^{k_1}}+\cdots+\frac{1}{2^{k_1+\cdots+k_r}}+\frac{1}{2^{k_1+\cdots+k_r+l-1}}        .$$
          Since 
          \[
          \begin{split}
          &\;\;\;\; F_n\left(\tau((k_1,k_2,\cdots,k_r,l, \{1\}^\infty   )     )\right)- T_n(k_1,k_2, \cdots,k_r)\\
          &=\mathop{\int}_{[0,1]^{k_1+\cdots+k_r+l}}   \left[1+t_n(x_1\cdots x_{k_1})      \right] \cdots  \left[1+t_n(x_1\cdots x_{k_1+\cdots+k_r})      \right] \Bigg{[} \left(1+t_n(x_1\cdots x_{k_1+\cdots+k_r+l})      \right)\\
          &\;\;\;\;\cdot \mathrm{e}_n\Big{(} x_1\cdots x_{k_1+\cdots+k_r+l}    \Big{)}    -1   \Bigg{]}  dx_1dx_2\cdots dx_{k_1+\cdots+k_r+l} \\
           &<\mathop{\int}_{[0,1]^{k_1+\cdots+k_r+l}}   \left[1+nx_1\cdots x_{k_1}      \right] \cdots  \left[1+n x_1\cdots x_{k_1+\cdots+k_r}      \right] \Bigg{[} \left(1+n x_1\cdots x_{k_1+\cdots+k_r+l}      \right)\\
          &\;\;\;\;\cdot \mathrm{exp}\,\Big{(} \left( 1+\mathrm{log}\,n   \right)x_1\cdots x_{k_1+\cdots+k_r+l}  \Big{)}    -1   \Bigg{]}  dx_1dx_2\cdots dx_{k_1+\cdots+k_r+l} \\
        &<\mathop{\int}_{[0,1]^{k_1+\cdots+k_r+l}}   \left[1+nx_1\cdots x_{k_1}      \right] \cdots  \left[1+n x_1\cdots x_{k_1+\cdots+k_r}      \right] \Bigg{[} \left(1+n x_1\cdots x_{k_1+\cdots+k_r+l}      \right)\\
          &\;\;\;\;\cdot \Big{(}  1+\sum_{m\geq 1}\frac{(1+\mathrm{log}\,n)^m}{m!}   x_1\cdots x_{k_1+\cdots+k_r+l}  \Big{)}    -1   \Bigg{]}  dx_1dx_2\cdots dx_{k_1+\cdots+k_r+l} \\    \end{split}
             \]
             \[
             \begin{split}  
             &<\left( n+ \sum_{m\geq 1} \frac{(1+\mathrm{log}\,n)^m}{m!} + n\cdot \sum_{m\geq 1} \frac{(1+\mathrm{log}\,n)^m}{m!}  \right)\mathop{\int}_{[0,1]^{k_1+\cdots+k_r+l}}   \left[1+nx_1\cdots x_{k_1}      \right] \cdots \\
             &\;\;\;\; \cdot \left[1+n x_1\cdots x_{k_1+\cdots+k_r}      \right]  x_1\cdots x_{k_1+\cdots+k_r+l}   \cdot  dx_1dx_2\cdots dx_{k_1+\cdots+k_r+l} \\   
                              &<\left( n+ \sum_{m\geq 1} \frac{(1+\mathrm{log}\,n)^m}{m!} + n\cdot \sum_{m\geq 1} \frac{(1+\mathrm{log}\,n)^m}{m!}  \right) \frac{e^{2n}}{2^{k_1+\cdots+k_r+l}}          \end{split}
          \] 
          and 
          \[
          \begin{split}
&\;\;\;\; H_n(\tau((k_1,k_2,\cdots,k_r,l, \{1\}^\infty   )     ) ) -  H_n(x_0 )-h_n(x_0)      \\
&=       \sum_{\substack{  x_0<y<w_l  \\   y\in \tau(\widehat{\mathcal{T}}^R)   }} h_n(y)   \\
&=\sum_{q_1\geq l}\Bigg{[} h_n\left(  \frac{1}{2^{k_1}}+\cdots+\frac{1}{2^{k_1+\cdots+k_r}}+\frac{1}{2^{k_1+\cdots+k_r+q_1}}   \right)+  \sum_{m\geq 2}\sum_{q_2,\cdots,q_m\geq 1}   \\
&\;\;\;\;\;h_n\left(  \frac{1}{2^{k_1}}+\cdots+\frac{1}{2^{k_1+\cdots+k_r}}+\frac{1}{2^{k_1+\cdots+k_r+q_1}} +\cdots+  \frac{1}{2^{k_1+\cdots+k_r+q_1+\cdots+q_m}}     \right)   \Bigg{]}\\
&\leq 2n e^{(n+2)n} \sum_{q_1\geq l}\Bigg{(}   \frac{1}{(n+2)^{k_1+\cdots+k_r+q_1}}  +  \sum_{m\geq 2}\sum_{q_2,\cdots,q_m\geq 1}   \frac{1}{(n+2)^{k_1+\cdots+k_r+q_1+\cdots+q_m}}    \Bigg{)}\\
&\leq \frac{ 2n e^{(n+2)n}}{(n+2)^{k_1+\cdots+k_r} }\sum_{q_1\geq l} \frac{1}{ (n+2)^{q_1}     }\left(1+\sum_{m\geq 2}\frac{1}{(n+1)^{m-1}}    \right)\\
&\leq \frac{ 2n e^{(n+2)n}}{(n+2)^{k_1+\cdots+k_r} }  \cdot   \frac{1}{(n+2)^{l-1} (n+1)}  \cdot    \frac{n+1}{n}        \\
&\leq \frac{2 e^{(n+2)n}}{(n+2)^{k_1+\cdots+k_r+l-1}},\\
\end{split}
\]      
one has 
\[
\begin{split}
&\;\;\;\;\mathop{\lim}_{m\rightarrow +\infty} G_n(u_m)\\
&=T_n(k_1,k_2, \cdots,k_r)+H_n(x_0 )+h_n(x_0)\\
&= \eta_n\left((k_1,\cdots,k_{r-1}, k_r+1, \{1\}^\infty         )       \right)    +H_n(x_0)\\
&=   F(x_0)+H_n(x_0)\\
&= G_n(x_0).  
\end{split}
\]

In conclusion, $G_n$ is continuous on  $(0,1]$. The above analysis also reveals that the function $H_n$ is left continuous on $(0,1)$.

Since $G_n$ is continuous and $\tau(\widehat{\mathcal{T}}^R)$ is dense in $(0,1]$, to show that $G_n$ is strictly increasing, it suffices to show that $G_n$ is strictly increasing on $\tau(\widehat{\mathcal{T}}^R)$.

By the observation in Remark \ref{sleq}, it is easy to check that 
\[
\eta_n((k_1,\cdots,k_r,\cdots))> \eta_n((l_1,\cdots,l_r,\cdots))\]
if $k_i\leq l_i$ for all $i\geq 1$ and $ (k_1,\cdots,k_r,\cdots)\neq (l_1,\cdots,l_r,\cdots)$.

For $r,l\geq 1, k_1,\cdots, k_r,k_{r+1}\geq 1$, define 
\[
y_1=\frac{1}{2^{k_1}}+\cdots+\frac{1}{ 2^{k_1+\cdots+k_r}}+\frac{1}{ 2^{k_1+\cdots+k_r+k_{r+1}}},
\]
\[
y_2=\frac{1}{2^{k_1}}+\cdots+\frac{1}{ 2^{k_1+\cdots+k_r}},\]
\[
y^{(l)}_3=\frac{1}{2^{k_1}}+\cdots+\frac{1}{ 2^{k_1+\cdots+k_{r-1}}}+\frac{1}{ 2^{k_1+\cdots+k_{r-1}+k_r+1}}+\cdots+\frac{1}{ 2^{k_1+\cdots+k_{r-1}+k_r+l}} .\]

Then we have \[
\begin{split}
&\;\;\;\; F_n(y_1)-F_n(y_2)\\
&= \eta_n\circ \tau\left((k_1,\cdots,k_r,k_{r+1}+1,\{1\}^\infty)\right)-  \eta_n\circ \tau\left((k_1,\cdots,k_{r-1},k_{r}+1,\{1\}^\infty)\right)     \\
&=\mathop{\int}_{[0,1]^{k_1+\cdots+k_r+k_{r+1}} } \left[1+t_n(x_1\cdots x_{k_1})\right]\cdots \left[ 1+t_n(x_1\cdots x_{k_1+\cdots+k_r})\right]
  e_n\left( x_1\cdots x_{k_1+\cdots+k_{r+1}}   \right)\\
  &\;\;\;\;\cdot dx_1dx_2\cdots dx_{k_1+\cdots+k_{r+1}}\\
&\;\;\;\;-  \mathop{\int}_{[0,1]^{k_1+\cdots+k_r} } \left[1+t_n(x_1\cdots x_{k_1})\right]\cdots \left[ 1+t_n(x_1\cdots x_{k_1+\cdots+k_{r-1}})\right]e_n\left( x_1\cdots x_{k_1+\cdots+k_{r}}   \right)\\
&\;\;\;\;\cdot  dx_1dx_2\cdots dx_{k_1+\cdots+k_{r}}\\
&= \mathop{\int}_{[0,1]^{k_1+\cdots+k_r+k_{r+1}} } \left[1+t_n(x_1\cdots x_{k_1})\right]\cdots \left[ 1+t_n(x_1\cdots x_{k_1+\cdots+k_{r-1}})\right]        \\
& \bigg{\{} \left[ 1+t_n(x_1\cdots x_{k_1+\cdots+k_r})     \right]   e_n\left( x_1\cdots x_{k_1+\cdots+k_{r+1}}   \right)-  e_n\left( x_1\cdots x_{k_1+\cdots+k_{r}}   \right)  \bigg{\}} dx_1\cdots dx_{k_1+\cdots+k_{r+1}}\\
\end{split}
\]
and 
\[
\begin{split}
&\;\;\;\; H_n(y_1)-H_n(y_2)\\
&=h_n(y_2)+  \sum_{\substack{y_2<y<y_1     \\ y\in\tau(\widehat{T}^R)}} h_n(y)\\
&= \mathop{\int}_{[0,1]^{k_1+\cdots+k_r}}\left[ 1+t_n(x_1\cdots x_{k_1})\right] \cdots \left[   1+t_n(x_1\cdots x_{k_1+\cdots +k_{r-1}})  \right]\Delta_n(x_1\cdots x_{k_1+\cdots+k_r})\\
    &\;\;\;\;\cdot dx_1\cdots dx_{k_1+\cdots+k_r}+\sum_{\substack{y_2<y<y_1     \\ y\in\tau(\widehat{T}^R)}} h_n(y) .\\
    \end{split}
\]
As a result, one has 
\[
\begin{split}
&\;\;\;\; G_n(y_1)-G_n(y_2)\\
&= F_n(y_1)+H_n(y_1)-F_n(y_2)-H_n(y_2)       \\
&=  \sum_{\substack{y_2<y<y_1     \\ y\in\tau(\widehat{T}^R)}} h_n(y)+   \mathop{\int}_{[0,1]^{k_1+\cdots+k_r+k_{r+1}} } \left[1+t_n(x_1\cdots x_{k_1})\right]\cdots \left[ 1+t_n(x_1\cdots x_{k_1+\cdots+k_{r-1}})\right]        \\
& \bigg{\{} \left[ 1+t_n(x_1\cdots x_{k_1+\cdots+k_r})     \right] \left[  e_n\left( x_1\cdots x_{k_1+\cdots+k_{r+1}}   \right)-1\right]    \bigg{\}} dx_1\cdots dx_{k_1+\cdots+k_{r+1}}   \\
&>  \sum_{\substack{y_2<y<y_1     \\ y\in\tau(\widehat{T}^R)}} h_n(y) + \mathop{\int}_{[0,1]^{k_1+\cdots+k_r+k_{r+1}} }  \left[e_n\left( x_1\cdots x_{k_1+\cdots+k_{r+1}}   \right)-1\right] dx_1\cdots dx_{k_1+\cdots+k_{r+1}}  \\
&>0.
\end{split}
\]
Similarly, one has 
\[
\begin{split}
&\;\;\;\; G_n(y_2)-G_n(y_3^{(l)})\\
&= F_n(y_2)-F_n(y_3^{(l)})+H_n(y_2)-H_n(y_3^{(l)})       \\
&=\mathop{\int}_{[0,1]^{k_1+\cdots+k_r} } \left[1+t_n(x_1\cdots x_{k_1})\right]\cdots \left[ 1+t_n(x_1\cdots x_{k_1+\cdots+k_{r-1}})\right]  e_n\left( x_1\cdots x_{k_1+\cdots+k_{r-1}+k_r}   \right)  \\
&\;\;\;\;\cdot dx_1\cdots dx_{k_1+\cdots +k_r}-  \mathop{\int}_{[0,1]^{k_1+\cdots+k_r+l} } \left[1+t_n(x_1\cdots x_{k_1})\right]\cdots \left[ 1+t_n(x_1\cdots x_{k_1+\cdots+k_{r-1}})\right]  \\
&\cdot \left[ 1+t_n(x_1\cdots x_{k_1+\cdots+k_{r}+1})\right] \cdots  \left[ 1+t_n(x_1\cdots x_{k_1+\cdots+k_{r}+l-1})\right]   e_n\left( x_1\cdots x_{k_1+\cdots+k_{r}+l}   \right)    \\
&\;\;\;\;\cdot dx_1\cdots dx_{k_1+\cdots +k_r+l}+h_n(y_3^{(l)}) + \sum_{\substack{y_3^{(l)}<y<y_2     \\ y\in\tau(\widehat{T}^R)}} h_n(y)\\
&>\mathop{\int}_{[0,1]^{k_1+\cdots+k_r} } \left[1+t_n(x_1\cdots x_{k_1})\right]\cdots \left[ 1+t_n(x_1\cdots x_{k_1+\cdots+k_{r-1}})\right]  e_n\left( x_1\cdots x_{k_1+\cdots+k_{r-1}+k_r}   \right)  \\
&\;\;\;\;\cdot dx_1\cdots dx_{k_1+\cdots +k_r}-  \mathop{\int}_{[0,1]^{k_1+\cdots+k_r+l} } \left[1+t_n(x_1\cdots x_{k_1})\right]\cdots \left[ 1+t_n(x_1\cdots x_{k_1+\cdots+k_{r-1}})\right]  \\
&\;\;\;\;\cdot \left[ 1+t_n(x_1\cdots x_{k_1+\cdots+k_{r}+1})\right] \cdots  \left[ 1+t_n(x_1\cdots x_{k_1+\cdots+k_{r}+l})\right]    dx_1\cdots dx_{k_1+\cdots +k_r+l}\\
&\geq \mathop{\int}_{[0,1]^{k_1+\cdots+k_r+l} } \left[1+t_n(x_1\cdots x_{k_1})\right]\cdots \left[ 1+t_n(x_1\cdots x_{k_1+\cdots+k_{r-1}})\right]  \Big{[} e_n\left( x_1\cdots x_{k_1+\cdots+k_{r-1}+k_r}   \right)\\
&- \left[ 1+t_n(x_1\cdots x_{k_1+\cdots+k_{r}+1})\right] \cdots  \left[ 1+t_n(x_1\cdots x_{k_1+\cdots+k_{r}+l})\right]  \Big{]}  dx_1\cdots dx_{k_1+\cdots +k_r+l}\\
&>0. \end{split}
\]
There the last inequality follows from Lemma \ref{111}.
In general cases, for  $x_1>x_2$  and $x_1,x_2\in\tau( \widehat{\mathcal{T}}^R)$, there is an $r_0\geq 1$ such that
\[
k_{r_0}<j_{r_0},
\]
\[
x_1=\frac{1}{2^{k_1}}+\cdots+\frac{1}{2^{k_1+\cdots+k_{r_0-1}}}+ \frac{1}{2^{k_1+\cdots+k_{r_0-1}+k_{r_0}}}+\cdots+\frac{1}{2^{k_1+\cdots+k_l}},
\]
\[
x_2= \frac{1}{2^{k_1}}+\cdots+\frac{1}{2^{k_1+\cdots+k_{r_0-1}}}+ \frac{1}{2^{k_1+\cdots+k_{r_0-1}+j_{r_0}}}+\cdots+\frac{1}{2^{k_1+\cdots+j_{m}}}.
\]
By the above analysis and induction, for $l\geq 1$, one has 
\[
\begin{split}
&\;\;\;\;G_n(x_1)\\
&>G_n\left(  \frac{1}{2^{k_1}}+\cdots+\frac{1}{2^{k_1+\cdots+k_{r_0-1}}}+ \frac{1}{2^{k_1+\cdots+k_{r_0-1}+k_{r_0}}}  \right)\\
&> G_n\left(  \frac{1}{2^{k_1}}+\cdots+\frac{1}{2^{k_1+\cdots+k_{r_0-1}}}+ \frac{1}{2^{k_1+\cdots+k_{r_0-1}+j_{r_0}}}\left(1+\cdots   +  \frac{1}{2^{ l}}\right) \right)   .     \\
\end{split}
\]
As for $l\geq j_{r_0+1}+\cdots+j_m$, we have 
\[
G_n\left(  \frac{1}{2^{k_1}}+\cdots+\frac{1}{2^{k_1+\cdots+k_{r_0-1}}}+ \frac{1}{2^{k_1+\cdots+k_{r_0-1}+j_{r_0}}}+\cdots   +  \frac{1}{2^{k_1+\cdots+k_{r_0-1}+j_{r_0}+ l}} \right)  \geq G_n(x_2)\]
by definition. 
In a word, we have $G_n(x_1)>G_n(x_2)$  for  $x_1>x_2$  and $x_1,x_2\in\tau( \widehat{\mathcal{T}}^R)$. Thus from the fact that $G_n$ is continuous on $(0,1]$ one has 
\[
G_n(x_1)>G_n(x_2), \;\;\;\; x_1>x_2, x_1,x_2\in (0, 1].
\]

$\hfill\Box$\\ 
   
  Now  we are ready to prove that \[\eta_n(\widehat{\mathcal{T}}- \widehat{\mathcal{T}}^{R})=(1,e^{1+\frac{1}{2}+\cdots+\frac{1}{n}}].\]
  For $y_0\in (1, e^{1+\frac{1}{2}+\cdots+\frac{1}{n}})$,  define 
  \[
  x_0=\mathop{\mathrm{inf}}\big{\{} x\in (0,1)\,\big{|}\, F_n(x)>y_0\big{\}}.
  \]
  By definition, for $\epsilon>0$, one has 
  \[
  F_n(x_0-\epsilon)=G_n(x_0-\epsilon)-H_n(x-\epsilon)\leq y_0
  \]
  Since $G_n$ is continuous and $H_n$ is left continuous, we have 
  \[
  F_n(x_0)\leq y_0.
  \]
  If $F_n(x_0)<y_0$, then 
  \[
  G_n(x_0)<H_n(x_0)+y_0.
  \]
  As $G_n$ is continuous, there is an $\epsilon >0$ such that 
  \[
  G_n(y)< H_n(x_0)+y_0
  \] 
  for all $y\in (x_0, x_0+\epsilon)$. Thus 
   \[
  G_n(y)< H_n(y)+y_0
  \] 
  for all $y\in (x_0, x_0+\epsilon)$.   This contradicts the definition of $x_0$.  So we have $F_n(x_0)=y_0$
  and 
  \[
  \eta_n(\widehat{\mathcal{T}})=(1,e^{1+\frac{1}{2}+\cdots+\frac{1}{n}}].  \]
  
  Furthermore,  for $y_0\in (1, e^{1+\frac{1}{2}+\cdots+\frac{1}{n}})$ define 
  \[
  x_1=\mathop{\mathrm{sup}}\big{\{} x\in (0,1)\,\big{|}\, F_n(x)=G_n(x)-H_n(x)=y_0\big{\}}.  \]
  Since $F_n$ is left continuous on $(0,1)$, one has $F_n(x_1)=y_0$.
  Assuming that $x_1 \in \tau(\widehat{\mathcal{T}}^{R})$,  from the fact that:\\
  $(i)$  $H_n$ is strictly increasing and discontinuous at $x_1$,\\
  $(ii)$  $G_n$ is continuous,\\
   there is an $\epsilon >0$ such that for $y\in (x_1,x_1+\epsilon)$,
  \[
  F_n(y)=G_n(y)-H_n(y)<y_0 . \]
  Define 
  \[
  x_2=\mathop{\mathrm{inf}}\big{\{} x\in (x_1+\frac{\epsilon}{2},1)\,\big{|}\, F_n(x)>y_0\big{\}}.
  \]
  By the similar trick we have $F_n(x_2)=y_0$ and $x_2>x_1+\frac{\epsilon}{2}$, this contradicts the definition of $x_1$. So we have $x_1 \in \tau(\widehat{\mathcal{T}}-\widehat{\mathcal{T}}^{R})$   and     \[\eta_n(\widehat{\mathcal{T}}- \widehat{\mathcal{T}}^{R})=(1,e^{1+\frac{1}{2}+\cdots+\frac{1}{n}}].\]

        \section*{Acknowledgements}
        The author wants to thank Shengyou Wen and Yufeng Wu for helpful discussions about fractal geometry.  This project is  supported by  the Natural Science Foundation of Hunan Province, China (Grant No. 2023JJ40691)  and  the National Natural Science Foundation of China (Grant No. 12201642).

\end{document}